\documentclass[12pt]{article}
\usepackage[final]{epsfig}
\usepackage{graphics}
\usepackage{amsmath}
\usepackage{amsfonts}
\usepackage{latexsym}
\usepackage{amssymb}
\usepackage{graphicx}
\usepackage{url}

\newtheorem{lemma}{Lemma}[section]
\newtheorem{proposition}[lemma]{Proposition}
\newtheorem{remark}[lemma]{Remark}

\newtheorem{theorem}{Theorem}

\newtheorem{corollary}[lemma]{Corollary}

\newtheorem{conjecture}{Conjecture}

\newcommand{\g}{{\gamma}}
\newcommand{\G}{{\Gamma}}
\newcommand{\proofend}{$\Box$\bigskip}

\newcommand{\R}{{\mathbb R}}

\newcommand{\HH}{{\mathbb {H}}} 

\def\proof{\paragraph{Proof.}}

\begin{document}

\title{On the bicycle transformation and the filament equation: results and conjectures}

\author{Serge Tabachnikov\footnote{
Department of Mathematics,
Pennsylvania State University, 
University Park, PA 16802, USA, 
tabachni@math.psu.edu}
}

\date{}
\maketitle

\section{Bicycle model, bicycle monodromy, and bicycle correspondence} \label{transf}

This paper concerns a simple model of bicycle kinematics and continues the work done in \cite{FLT,LT,Le,Ta,TT}. In this model, a bicycle is represented by an oriented segment of constant length $\ell$ that can move in such a way that the velocity of its rear end is aligned with the segment (the rear wheel is fixed on the bicycle frame). The motion takes place in the Euclidean plane but, equally well, albeit less physically, we consider the motion in $\R^n$ (see \cite{HPZ} for bicycle motion in the hyperbolic and elliptic planes).

The configuration space of oriented segments of fixed length in $\R^n$ is the spherization of the tangent bundle $ST\R^n$, and the non-holonomic constraint that we have imposed on the motion defines an $n$-dimensional distribution ${\cal D}$ therein. The motion of the bicycle is a smooth curve in $ST\R^n$ tangent to the distribution ${\cal D}$. 

The two projections $ST\R^n \to \R^n$, to the front and to the rear ends of the segment, yield the front and rear bicycle tracks. The former projection is transverse to ${\cal D}$, and the front track is a smooth curve, whereas the direction of the latter projection is contained in ${\cal D}$, and the rear track may have singularities (generically, semi-cubical cusps). This happens when the front wheel of the bicycle turns $90^{\circ}$. We denote the rear track by $\g$ and the front track by $\G$.

In this section, we review the relevant results from \cite{LT,Ta,TT} that will be used later.

The bicycle segment gives the rear track $\g$ coorinetation (the orientation of the normal bundle), well defined even at cusps. Conversely, given a cooriented curve $\g$, possibly,  with cusps, one uniquely reconstructs the respective curve $\G$ as the locus of endpoints of the tangent segments of length $\ell$. That is, the rear track and the direction of the motion determine the front track. 

On the other hand, the front track $\G$ determines the motion of the bicycle, and hence the rear track $\g$, only after its initial position is chosen (this is obvious geometrically, and it follows from the fact that the motion of the bicycle with a given front track is described by a first order differential equation, see \cite{LT} and equation (\ref{bidiff}) below). Thus the {\it bicycle monodromy} map $M_{\G,\ell}$ arises that takes the initial position of the bicycle to the terminal one. 

Assume that $\G$ is a closed curve. Then $M_{\G,\ell}$ is a self-map of the sphere $S^{n-1}$, well defined up to conjugation (due to the freedom of choice of the initial point on $\G$). 
Consider $S^{n-1}$ as the sphere at infinity of the hyperbolic space $\HH^n$. Then the M\"obius group $O(n,1)$ of isometries of  $\HH^n$ acts on $S^{n-1}$. The next result was proved  in \cite{LT}. 

\begin{theorem} \label{Moeb}
The bicycle monodromy $M_{\G,\ell}$ is a M\"obius transformation.
\end{theorem}

\begin{figure}[hbtp]
\centering
\includegraphics[width=2.1in]{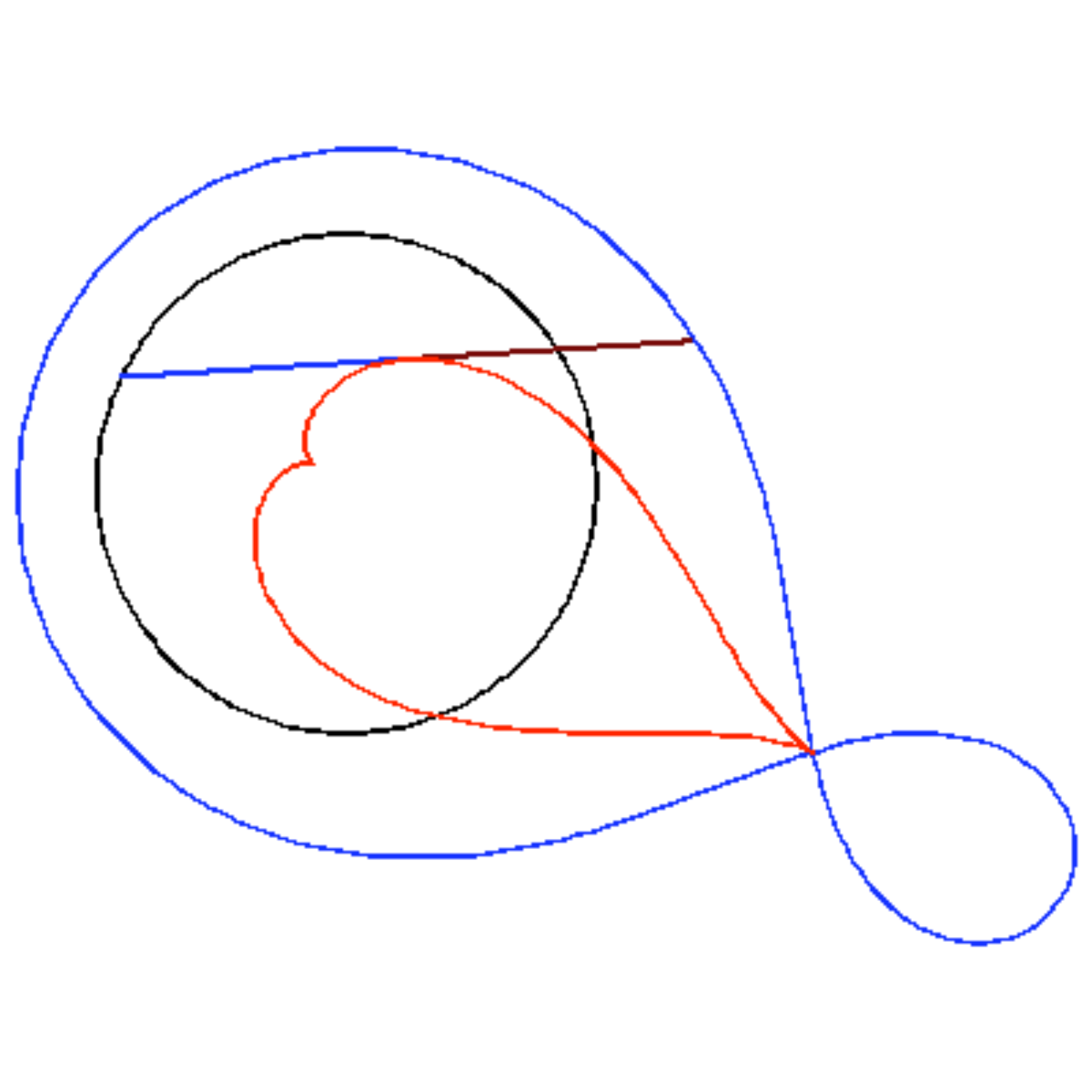}\quad
\includegraphics[width=2.4in]{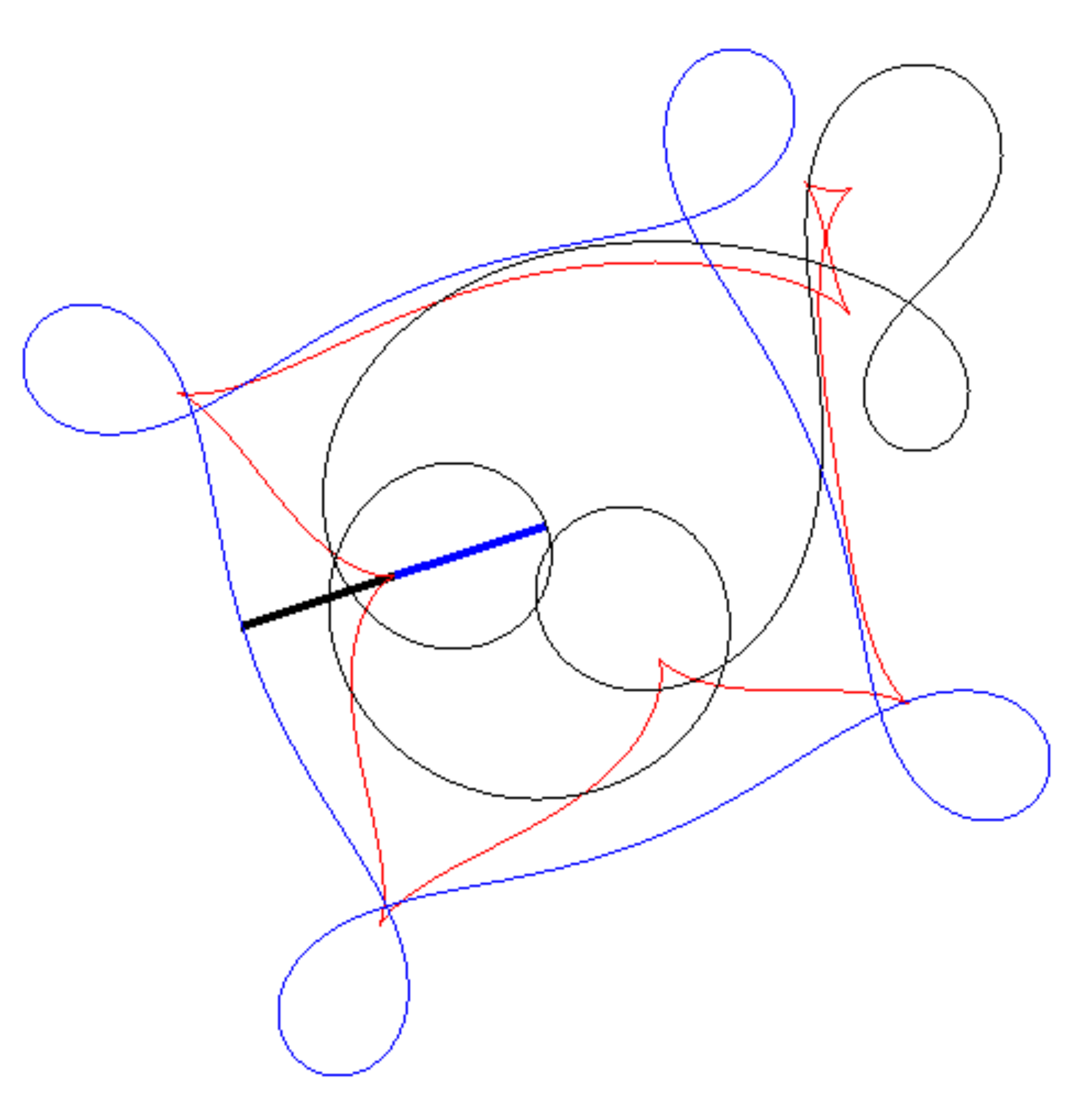}
\caption{Bicycle correspondence (figures courtesy of R. Perline).}
\label{Darboux}
\end{figure}

Let $\g$ be a closed rear track. It has two coorientations, and they correspond to two front tracks, $\G_1$ and $\G_2$ (one may think of a forward-facing and a backward-facing bicycle joined together). We say that $\G_1$ and $\G_2$ are in the {\it bicycle correspondence} and write ${\cal B}_{2\ell}(\G_1,\G_2)$. Equivalently, two points, $x_1$ and $x_2$, traverse the curves $\G_1$ and $\G_2$ in such a way that the distance $x_1x_2$ is equal to $2\ell$, and the velocity of the midpoint of the segment $x_1x_2$ is aligned with $x_1x_2$ (it follows that the speeds of $x_1$ and $x_2$ are equal). See Figure \ref{Darboux}.  

\begin{figure}[hbtp]
\centering
\includegraphics[width=1.9in]{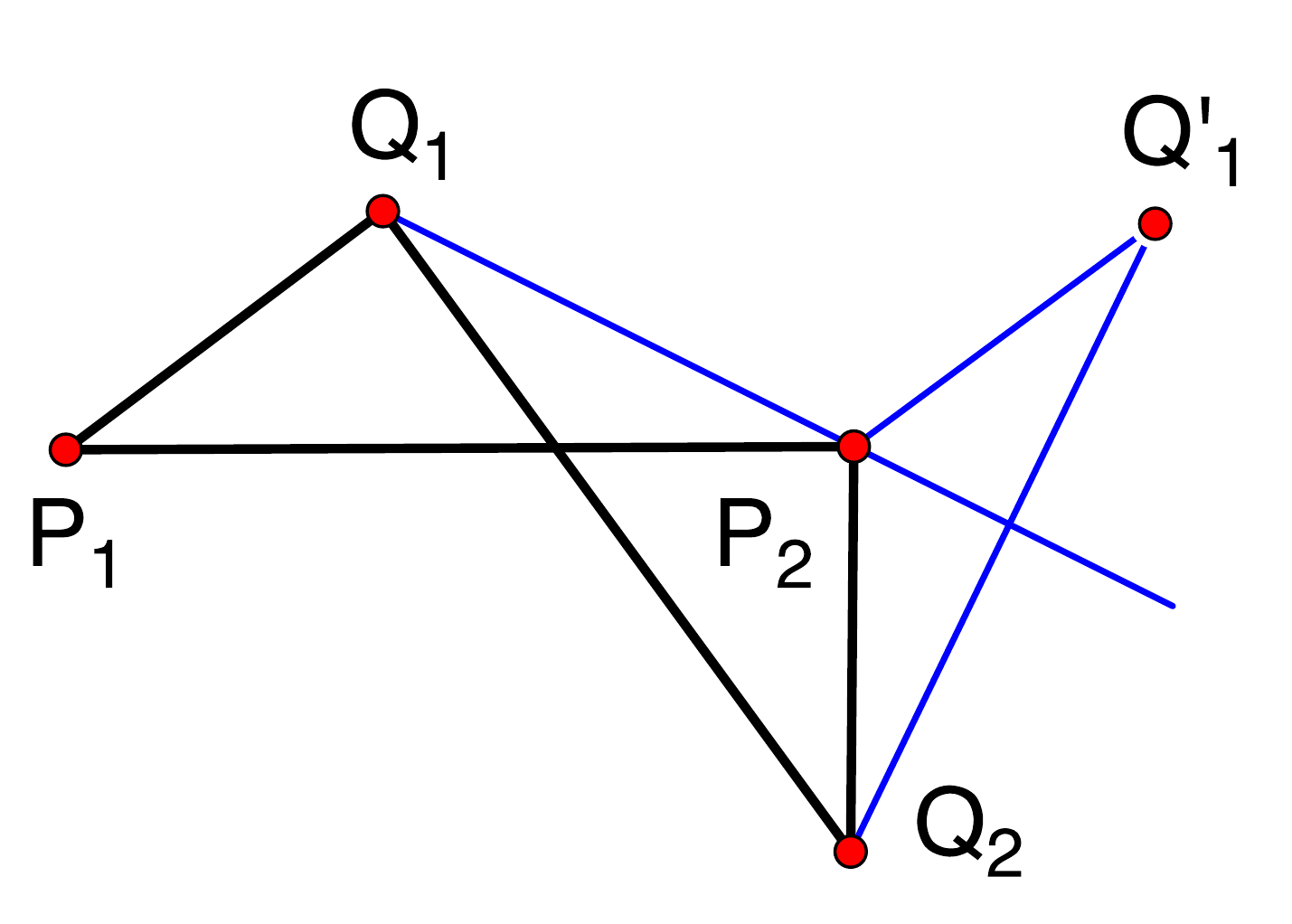}
\caption{Discrete bicycle transformation.}
\label{constr}
\end{figure}

The bicycle monodromy and the bicycle correspondence have discrete versions \cite{Ho,PSW,TT}. Let $P=(P_1,P_2,\ldots,P_n)$ be an $n$-gon in $\R^n$, and let $P_1 Q_1$ be a segment of length $2\ell$. Let $Q_2$ be the point such that $P_1 Q_1 P_2 Q_2$ is an isosceles trapezoid. In other words, one parallel translated the segment $P_1 Q_1$ to  $P_2 Q_1'$ and then reflects $Q_1'$ in the line $Q_1 P_2$, see Figure \ref{constr}. Continuing this, one constructs a polygon $Q_1,Q_2,\ldots,Q_n$.

\begin{figure}[hbtp]
\centering
\includegraphics[width=2.9in]{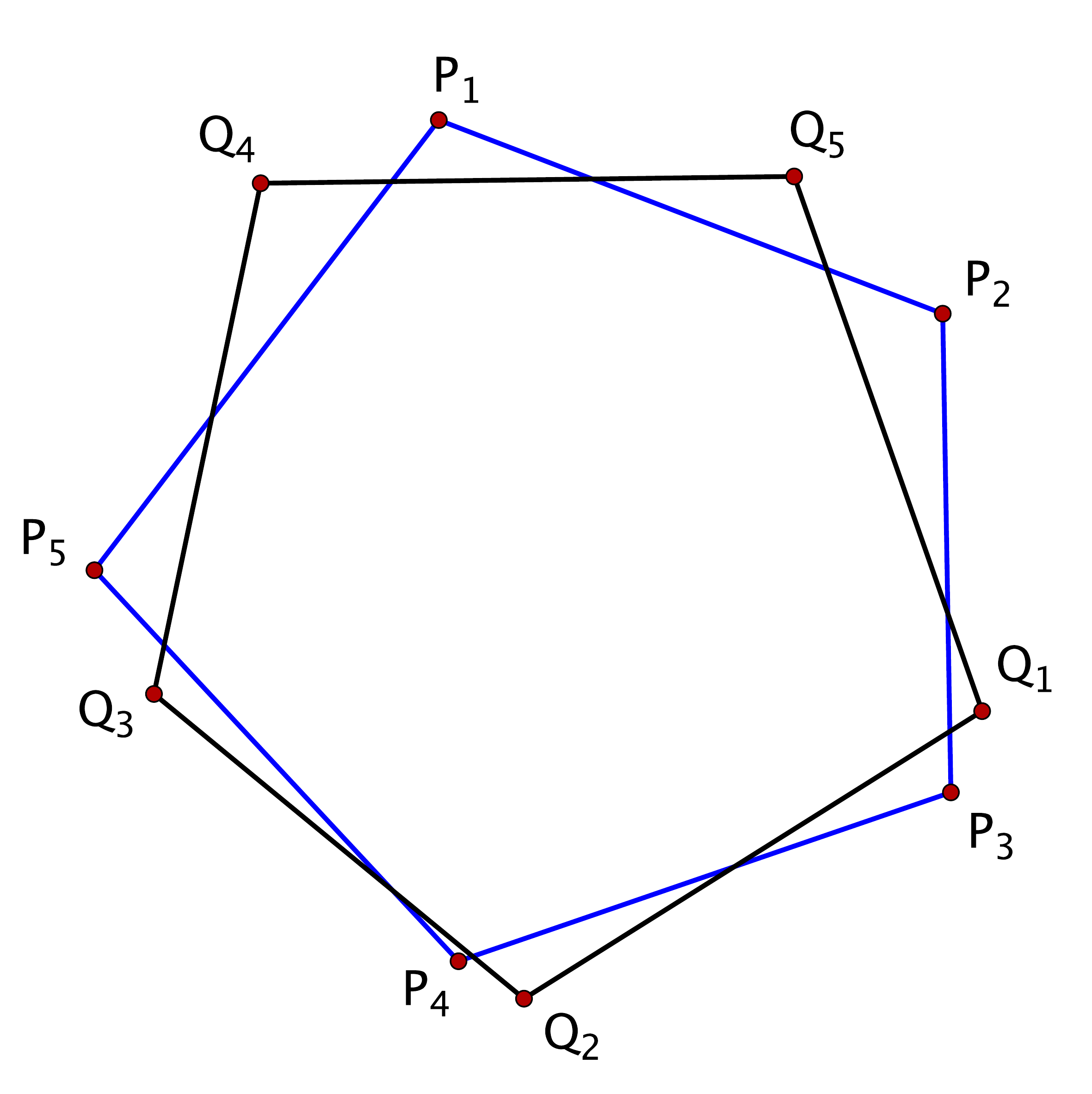}
\caption{Two polygons in the discrete bicycle correspondence.}
\label{chain}
\end{figure}

In the limit $n\to\infty$, the two polygons become curves in the bicycle correspondence. The point $Q_{n+1}$ does not necessarily coincide with $Q_1$, and we obtain the discrete bicycle monodromy $Q_1 \mapsto Q_{n+1}$. Theorem \ref{Moeb} holds in this case as well, see \cite{TT}. If the polygon $Q$ closes up, we say that the polygons $P$ and $Q$ are in the discrete bicycle correspondence, see Figure \ref{chain}.

The next two theorems hold both for closed curves and for polygons, see \cite{TT} for proofs in the discrete setting. We formulate the continuous versions.

\begin{theorem} \label{mono}
If closed curves $\G_1$ and $\G_2$ are in the bicycle correspondence, ${\cal B}_{\ell}(\G_1,\G_2)$, then the bicycle monodromies $M_{\G_1,\lambda}$ and $M_{\G_2,\lambda}$ are conjugated for every value of $\lambda$.
\end{theorem}

\begin{theorem} [Bianchi permutability] \label{permut}
Let $\G_1,\G_2$ and $\G_3$ be three closed curves such that ${\cal B}_{\ell}(\G_1,\G_2)$ and ${\cal B}_{\lambda}(\G_1,\G_3)$ hold. Then there exists a closed curve $\G_4$ such that ${\cal B}_{\lambda}(\G_2,\G_4)$ and ${\cal B}_{\ell}(\G_3,\G_4)$ hold.
\end{theorem}

The bicycle correspondence is not a mapping: for example,  the bicycle monodromy of a curve $\G_1$ may be the identity, and then there is a whole sphere worth of corresponding curves $\G_2$.
For instance, this happens when $\G_1$ is a circle, see Figure \ref{Darboux} on the left.

According to Theorem \ref{Moeb}, in dimension three, the bicycle monodromy is a complex fractional-linear transformation and, for a generic curve and generic value of the legth, it has two fixed points. Therefore the bicycle correspondence is a 2-2 map and, up to a binary choice, one can consider it as a mapping of space curves. 

In the plane, the bicycle monodromy is a real fractional-linear transformation, and one has the usual trichotomy: it can be elliptic, parabolic, or hyperbolic. Let $\G_1$ be a plane curve such that $M_{\G_1,\ell}$ is hyperbolic, that is, it has two fixed points. Then one can choose one of the two respective curves, $\G_2$ such that ${\cal B}_{2\ell}(\G_1,\G_2)$.   By Theorem \ref{mono}, $M_{\G_2,\ell}$ is again hyperbolic, and there is a unique curve $\G_3$, different from $\G_1$, with ${\cal B}_{2\ell}(\G_2,\G_3)$, and so on. In this way the bicycle correspondence becomes a mapping of the space of plane curves.

Our ultimate goal is   the following conjecture which we formulate in somewhat loose terms.

\begin{conjecture} \label{mainconj}
The bicycle correspondence in $\R^n$ is Liouville integrable.
\end{conjecture}

This   means that the correspondence has an invariant Poisson bracket and 
an infinite collection of independent integrals in involution that is, in an appropriate sense, complete. 

Theorem \ref{mono} is a source of integrals. Let
$$
\det(E-t M_{\G,\lambda}) = \sum_{i=1}^{n+1} \sum_{j=1}^{\infty} I_{ij}(\G) t^i \lambda^j
$$
be the Taylor expansion of the characteristic polynomial of the monodromy, considered as an element of the M\"obius group $O(n,1)$. Then the functions $I_{ij}(\G)$ on the space of curves in $\R^n$ are conserved by the bicycle correspondence. We call them the {\it monodromy integrals}. 

\section{Two differential 2-forms} \label{forms}

Consider two differential 2-forms on the space of smooth embedded parameterized curves, the first defined in $\R^n$ for all $n$, and the second in $\R^3$:
$$
\omega (u,v) = \int u'(t)\cdot v(t)\ dt,\quad \Omega(u,v)=\int \det(\Gamma'(t),u(t),v(t))\ dt.
$$
Here $\Gamma(t)$ is a curve, and $u(t),v(t)$ are vector fields along $\Gamma$. 

Note that the integrals do not depend on the parameterization. Note also that $\Omega$ vanishes on the submanifold of plane curves. The form $\omega$ can be defined on any Riemannian manifold, and the form $\Omega$ depends only on the volume form (and not on the metric).

\begin{proposition} \label{closedforms}
The 2-forms $\omega$ and $\Omega$ are closed. The forms are degenerate: the 
kernel of $\omega$ consists of the parallel vector fields along the curve, and the kernel of $\Omega$ consists of the tangent vector fields.
\end{proposition}

\proof Consider the following 1-forms:
$$
\alpha(u)=\frac{1}{2} \int \Gamma'(t)\cdot u(t)\ dt,\ \beta(u)=-\frac{1}{3} \int \det(\Gamma(t),\Gamma'(t),u(t))\ dt.
$$
We claim that 
$$
\omega=d\alpha,\ \Omega=d\beta.
$$
Indeed, let $u$ and $v$ be commuting vector fields in a neighborhood of $\G$. One has 
$$
d\alpha(u,v)=L_u (\alpha(v)) - L_v (\alpha(u)) - \alpha([u,v]).
$$
The first two terms on the right hand side are easy to compute:
$$
L_u (\alpha(v)) = \int u'(t) \cdot v(t)\ dt,\  L_v (\alpha(u)) = \int v'(t) \cdot u(t)\ dt = - \int u'(t) \cdot v(t)\ dt,
$$
and the result for $\omega$ follows. A similar computation works for $\Omega$.

If $u$ is in the kernel of $\omega$, then $\int u'(t)\cdot v(t)\ dt =0$ for all vector fields $v$ along $\G$. It follows that $u'=0$, that is, $u$ is parallel along $\G$. 
Likewise, if $u$ is in the kernel of $\Omega$, then $\int \det(\Gamma'(t),u(t),v(t))\ dt =0$ for all vector fields $v$ along $\G$. This implies that $u$ is collinear with $\Gamma'$.
\proofend

\begin{corollary} \label{sympl}
The form $\omega$ descends on the quotient space of curves modulo parallel translations and is symplectic therein, and the form $\Omega$ descends on the quotient space of curves mod ${\rm Diff}_{+} (S^1)$, that is,  on the space of non-parameterized oriented curves, and  is a symplectic structure therein.
\end{corollary}

The latter is the celebrated Marsden-Weinstein symplectic structure, see   \cite{AK,Ca, MW}. One can identify isotopic equivalence classes of smooth knots with coadjoint orbits of the group of volume-preserving diffeomorphisms of $\R^3$ (the pairing of a curve with a divergence-free vector field is given by the flux across a film that spans the curve). Under this identification, the  Marsden-Weinstein symplectic structure coincides with the Kirillov-Kostant-Souriau symplectic structure on coadjoint orbits of a Lie group. 

Our main observation is as follows.

\begin{theorem} \label{bisymp}
The bicycle correspondence preserves the forms $\omega$ and $\Omega$, the former in any dimension, and the latter in dimension three.
\end{theorem}

\proof
Consider a smooth arc of the rear track $\gamma$ between two consecutive cusps and give it arc length parameterization with parameter $t\in[a,b]$. Let $\Gamma_{\pm}=\gamma\pm\ell\gamma'$ be the respective arcs of the two front  tracks which are in the bicycle correspondence.

Let $u(t),v(t)$ be two vector fields along $\gamma (t)$. The differentials of the maps $\gamma \mapsto \G_{\pm}$ take $u$ and $v$ to vector fields $U$ and $V$ along $\Gamma_{\pm}$. We claim that
\begin{equation} \label{uU}
U=u+\ell (u'- (\gamma'\cdot u') \gamma'),
\end{equation}
and likewise for $V$; here $\ell$ may be positive or negative, depending on whether one considers $\G_+$ or $\G_-$.

To prove (\ref{uU}), let $\gamma_1=\gamma+\varepsilon u$ be a variation of $\gamma$. Then $|\gamma_1'|=1+\varepsilon \gamma'\cdot u'$, and one has: 
$$
\Gamma_1=\gamma_1+\ell \frac{\gamma'_1}{|\gamma'_1|}=\gamma +\ell \gamma' + \varepsilon ( u + \ell (u'-(\gamma'\cdot u') \gamma')),
$$
as claimed (all computations are mod $\varepsilon^2$, and $\ell$ is signed). 

It follows from (\ref{uU}) that 
$$
U'=u'+\ell(u''-(\gamma''\cdot u')\gamma'-(\gamma'\cdot u'')\gamma'-(\gamma'\cdot u')\gamma'').
$$
We want to show that the part of the integral $\int U'\cdot V dt$ which is odd in $\ell$ vanishes. The only odd  term is linear in $\ell$; it is equal to
$$
\int_a^b [u'\cdot v'+u''\cdot v]- [(\gamma'\cdot v')(\gamma'\cdot u')+(\gamma'\cdot v)(\gamma''\cdot u')+(\gamma'\cdot v)(\gamma'\cdot u'')+(\gamma''\cdot v)(\gamma'\cdot u')]\ dt.
$$
The  integral of the first bracket is equal to $\ell (u'\cdot v)|_a^b$, and of the second to $\ell (\gamma'\cdot v)(\gamma'\cdot u')|_a^b$ 
(this would finish the proof if $\gamma$ had no cusps).

Now consider the next smooth arc of $\gamma$, arc length parameterized by $t\in[b,c]$. A similar computation applies to this arc, with $\ell$ changing the sign. The contribution of the cusp, corresponding to the value of parameter $b$, is $\ell (u'\cdot v)(b)$, computed for the first arc, and 
$\ell (u'\cdot v)(b)$, computed for the second one. They cancel each other:  the sign of $u'$ also changes to the opposite because the direction of motion reverses in the cusp. Likewise, one has cancelation of the terms $\ell (\gamma'\cdot v)(\gamma'\cdot u')$, and this implies that $\bar T_{\ell}^*(\omega)=\omega$.

Similarly, the computation for $\Omega$ involves
$
\int \det (\Gamma',U,V) dt.
$
There are two terms odd in $\ell$, the linear and the cubic ones. The linear  term is  $\int \det(\gamma',u,v)' dt$, again, the integral of a derivative. And again, the total integral of a full derivative vanishes, taking the sign changes at the cusps into account.

Finally, consider the cubic in $\ell$ term,
$$
\int \det(\gamma'',u'-(\gamma'\cdot u')\gamma',v'-(\gamma'\cdot v')\gamma')\ dt.
$$
The determinant is identically zero since all three vectors are orthogonal to  $\gamma'$, and hence linearly dependent. 
\proofend

The monodromy integrals descend to both quotient spaces, curves modulo parallel translations and non-parameterized curves.

\begin{conjecture} \label{commute}
The monodromy integrals Poisson commute on these quotient spaces with respect to the Poisson bracket induced, respectively, by the forms $\omega$ (in every dimension) and $\Omega$ (in dimension three). 
\end{conjecture}

\section{The bicycle transformation in dimension three and the filament equation} \label{compint}

The {\it filament (binormal, localized induction, smoke ring) equation} is an evolution on curves in $\R^3$
$$
\dot \Gamma = \Gamma' \times \Gamma''
$$ 
where $\Gamma(x)$ is arc length parameterized  curve, prime is $d/dx$, and dot is the time derivative. 

This system (equivalent to the non-linear Schr\"odinger equation via the Hasimoto transformation)  is completely integrable in the following sense, see \cite{LP,La}. It is a Hamiltonian system with respect to the symplectic structure induced by $\Omega$ on the space of arc length parameterized curves, the Hamiltonian function being the perimeter length of the curve.
It has a hierarchy of Poisson commuting integrals $F_1,F_2,\dots$ that starts with
\begin{equation} \label{integrals}
\int 1\ dx,\ \int \tau\ dx,\ \int \kappa^2\ dx,\ \int \kappa^2\tau\ dx, \int \left((\kappa')^2+\kappa^2 \tau^2 -\frac{1}{4} \kappa^4\right)\ dx, \dots
\end{equation}
where $\tau$ is the torsion, and $\kappa$ is the curvature. One also has a hierarchy of commuting Hamiltonian vector fields $X_0, X_1, X_2,\dots$ that starts with
\begin{equation} \label{fields}
-T,\ \kappa B, \frac{\kappa^2}{2} T + \kappa' N + \kappa \tau B, \kappa^2\tau T +(2\kappa'\tau + \kappa\tau')N+\left(\kappa\tau^2-\kappa''-\frac{\kappa^3}{2}\right) B, \dots
\end{equation}
where $T,N,B$ is the Frenet frame along $\Gamma$, and where $X_i$ is the Hamiltonian vector field of $F_i$. The vector fields $X_i$ satisfy the recurrence relation
\begin{equation} \label{recur}
T\times X_n = X_{n-1}'.
\end{equation}

In dimension three, the bicycle transformation is completely integrable, in the same sense as the filament equation.

\begin{theorem} \label{int}
The functions  (\ref{integrals}) are integrals of the bicycle transformation (independently of the length parameter $\ell$). The bicycle transformation commutes with the vector fields (\ref{fields}).
\end{theorem}

\proof We claim that
\begin{equation} \label{ident}
\omega (X_{n-1},\cdot)  = \Omega (X_n,\cdot) = dF_n.
\end{equation}
Indeed, the first equality follows from (\ref{recur}) and the definitions of the 2-forms, and the second from the fact that $X_n$ is the Hamiltonian vector field of the function $F_n$.

Assume that the bicycle transformation preserves the field $X_{n-1}$. By Theorem \ref{bisymp}, it  preserves $\omega$, and hence by (\ref{ident}), it preserves $dF_n$. 
Assume now that the bicycle transformation preserves the differential $dF_n$. By Theorem \ref{bisymp}, it  preserves $\Omega$, and hence the Hamiltonian vector field $X_n$. 

Since the bicycle transformation preserves $dF_n$, it changes $F_n$ by a constant: $F_n \mapsto F_n +c$. To see that this constant is zero, it suffices to consider a circle of sufficiently large radius. Such a circle is fixed by the bicycle transformation, hence $c=0$.

Thus one may argue inductively, starting with the function $F_1$, the perimeter length,  clearly preserved by the bicycle transformation.
\proofend

Theorem \ref{int} is hardly new: the bicycle transformation is the B\"acklund transformation for the filament equation, see, e.g., \cite{CI,RS}. One of the reasons  to use the term `bicycle' (rather than B\"acklund, or Darboux) transformation is to decouple it from the filament equation which, unlike the bicycle transformation, is specifically 3-dimensional. 

\begin{conjecture} \label{depend}
In dimensions three and two, the monodromy integrals and the integrals (\ref{integrals}) are functionally dependent of each other.  
\end{conjecture}

Specifically, we conjecture that the integrals (\ref{integrals}) are the consecutive terms in the Taylor expansion at $\lambda =0$ of the conjugacy invariant 
$$
\frac{{\rm Tr}^2 (M_{\G,\lambda})}{\det (M_{\G,\lambda})}
$$ 
where the monodromy is considered as a $2\times 2$ complex (in dimension 3) or real (in dimension 2) matrix.

\begin{remark}
{\rm In dimension two, the motion of the bicycle is described by a differential equation
\begin{equation} \label{bidiff}
\frac{d\alpha(x)}{dx} + \frac{\sin\alpha(x)}{\ell}= \kappa(x),
\end{equation}
where $\alpha$ is the steering angle, $\kappa$ is the curvature of the front track $\G$, $x$ is the arc length parameter on the front track $\G$, and $\ell$ the length of the bicycle segment, see \cite{LT}. It is proved in \cite{LT} that the derivatives $\sigma_{1,2}$ of the monodromy at its fixed points are equal to $\exp (\pm \int \cos\alpha\ dx)$, where $\alpha(x)$ is a periodic solution of (\ref{bidiff}).
 One has
$$
\frac{1}{\sigma_1} + \frac{1}{\sigma_2} = \frac{{\rm Tr}^2 (M_{\G,\lambda})}{\det (M_{\G,\lambda})} -2,
$$
hence the quantity $\int \cos\alpha\ dx$ is an integral of the bicycle transformation. 

This suggest a method of obtaining integrals: consider a periodic solution $\alpha(x)$ of equation (\ref{bidiff}) as a function of $\ell$ and expand $\cos \alpha$ as a Taylor series in $\ell$. Then the integral of each coefficient of this series is a conserved quantity of the bicycle transformation. The first three integrals obtained this way are 
$$
\int 1\ dx, \int \kappa^2\ dx,\  \int \left((\kappa')^2 -\frac{1}{4} \kappa^4\right)\ dx,
$$
as predicted by Conjecture \ref{depend}.
}
\end{remark}

Now let $M$ be an oriented 3-dimensional Riemannian manifold. One can define the differential 2-forms $\omega$ and $\Omega$ and considers the filament equation on $M$. One can also define an analog of the bicycle correspondence: the bicycle frame is a geodesic segment of fixed length.
It is natural to ask when one has an analog of the recursion scheme (\ref{ident}) for the filament equation and  an analog of Theorem \ref{int}.

\begin{conjecture} \label{compatible}
This holds if and only if $M$ has constant curvature.
\end{conjecture}

\section{In the plane. Zindler curves} \label{plane}

A plane analog of the filament equation is the {\it planar filament equation}
$$
\dot \G = \frac{\kappa^2}{2} T + \kappa' N
$$
where $(T, N)$ is the planar Frenet frame. This system is equivalent to the modified Korteweg-deVries equation, the same way as the filament equation to the non-linear Schr\"odinger equation.

It is shown in \cite{LP1} that the planar filament equation has an infinite hierarchy of integrals, namely,  the odd-numbered ones in the sequence (\ref{integrals}) (the torsion $\tau$ vanishes for plane curves). Likewise, the odd-numbered vector fields in the list (\ref{fields}) commute and preserve the planarity of curves. With this modification, Theorem \ref{int} holds for  the bicycle transformation in dimension two.

\begin{figure}[hbtp]
\centering
\includegraphics[width=1.2in]{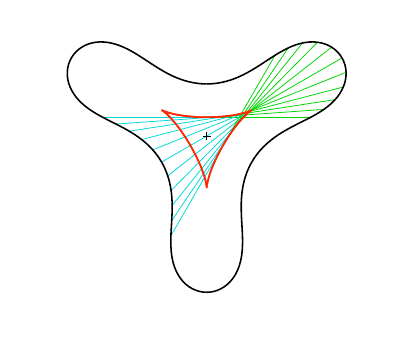}\ 
\includegraphics[width=1.2in]{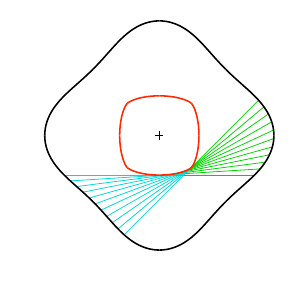}\ 
\includegraphics[width=1.2in]{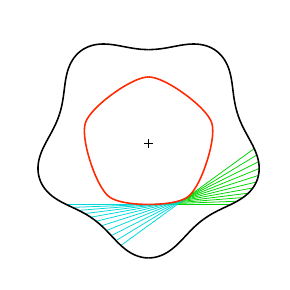}\ 
\includegraphics[width=1.2in]{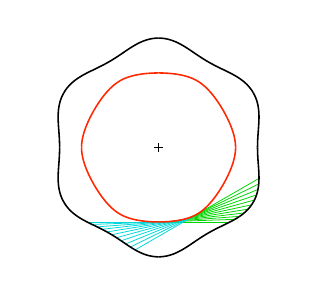}
\caption{Examples of Zindler curves from \cite{We6}.}
\label{Zindler}
\end{figure}

An interesting `bicycle' problem is to describe the situations when, given closed rear and  front tracks, $\g$ and $\G$, one cannot tell which way the bicycle went \cite{Fi}. A trivial example  is a pair of concentric circles; for numerous non-trivial examples, see \cite{We1}--\cite{We6} and Figure \ref{Zindler}.

Such an `ambiguous' front track $\G$ can be characterized by the property that a chord of fixed length $\ell$ can traverse the curve in such a way that the velocity of its mid-point is aligned with the chord. That is, $\G$ is in the bicycle correspondence with itself: ${\cal B}_{\ell}(\G,\G)$. Following \cite{BMO1,BMO2}, we call such curves {\it Zindler curves} (see \cite{Zi}).\footnote{Incidentally, (the interiors of) Zindler curves are  solutions to the 2-dimensional case of Ulam's problem: which bodies float in equilibrium in all positions? See \cite{Au}.} Of course, Zindler curves can be defined in other dimensions as well.

The 2-dimensional version of Theorem \ref{int} has the following corollary. 

\begin{corollary} \label{Zflow}
Let $\G$ be a Zindler curve, and let $\G_1$ be the result of evolving $\G$ along the planar filament flow. Then $\G_1$ is also a Zindler curve with the same length parameter $\ell$. The same is true 
if the planar filament field is replaced by  any of the odd-numbered vector fields in (\ref{fields}) or their linear combination.
\end{corollary}

\proof
As we mentioned earlier, the end points of the segment of fixed length move along a Zindler curve with the same speed. Let  $\G(x)$ be an arc length parameterization. Then ${\cal B}_{\ell} (\G(x),\G(x+c))$ for some constant $c$. Since the bicycle transformation commutes with the planar filament flow, and the latter preserves the arc length, we have ${\cal B}_{\ell} (\G_1(x),\G_1(x+c))$, as needed.
\proofend

Recall that a soliton of the filament equation is a curve that is a critical point for a linear combination of integrals of motion. 

\begin{conjecture} \label{soliton}
The Zindler curves are solitons of the planar filament equation.
\end{conjecture}

In particular, according to preliminary computations of R. Perline, the Zindrel curves constructed by Wegner \cite{We5} are pressurized elastica, that is, extrema of a linear combination of area, length, and total squared curvature. 

\section{Other integrals} \label{otherint}

The bicycle correspondence in $\R^n$  has additional integrals.

\begin{theorem} \label{other}
The bivector and the vector,
$$
A(\Gamma)=\int_{\Gamma} \Gamma(t)\wedge\Gamma'(t)\ dt\ \ {\rm and}\ \
J(\Gamma)=\int_{\Gamma} (\Gamma(t)\cdot\Gamma'(t))\ \Gamma(t)\ dt,
$$
are preserved by  the bicycle correspondence.
\end{theorem}

\proof As before,  $\Gamma=\gamma+\ell \gamma'$ with $t$ being an arc length parameter along $\gamma$. The odd in $\ell$ part of $A(\Gamma)$ is $\ell \int \gamma\wedge\gamma''\ dt$. Note that $\gamma\wedge\gamma'' = (\gamma\wedge\gamma')'$. This would prove that the integral vanishes for a smooth rear track $\gamma$. If $\gamma$ has cusps, we note that $\gamma'$ and $\ell$ change  sign in each cusp, and we conclude  similarly to the proof of Theorem \ref{bisymp}.

Likewise, for $J(\G)$, the odd in $\ell$ parts are
$$
\ell \int ((1+\gamma\cdot\gamma'')\ \gamma + (\gamma\cdot\gamma')\ \gamma')\ dt\ \ {\rm and}\ \ 
\ell^3 \int (\g' \cdot \g'')\ \g'\ dt.
$$
Concerning the first integral, the integrand is $((\gamma\cdot\gamma')\ \gamma)'$, and that $(\gamma\cdot\gamma')\ \gamma$ changes sign at cusps, along with $\ell$. The second integral vanishes because $\g(t)$ is an arc length parameterization, and $\g' \cdot \g''=0$.
\proofend  

If the ambient dimension   is $n$ then the area bivector $A(\Gamma)$ provides $n\choose 2$ integrals, the areas of the projections on all the coordinate 2-planes, and the ``centroid" vector $J(\G)$ provides $n$ integrals. 
In dimension two, the meaning of these integrals is as follows: $A(\G)$ is twice the area bounded by the curve $\G$, and the vector $J(\G)$, rotated $90^{\circ}$ and divided by the area, is the center of mass of the domain bounded by $\G$. 

We finish with a remark about the area vector as an integral of the filament equation.
Let $a$ be a constant vector field in $\R^3$ (that is in the kernel of the form $\omega$), and let 
$$
F_a = \frac{1}{2}A(\G) \cdot a = \frac{1}{2} \int_{\G} (\G(t) \times \G'(t)) \cdot a\ dt,
$$
the projection of the area vector along $a$. Likewise, define
$$
J_a = - J(\G) \cdot a = -  \int_{\G} (\Gamma(t)\cdot\Gamma'(t))\ (\Gamma(t)\cdot a)\ dt.
$$

\begin{lemma} \label{proj}
One has: 
$$
i_a \Omega= dF_a = i_{\G\times a} \omega,\ i_{\G\times a} \Omega = dJ_a.
$$
\end{lemma}

\proof Let $u$ be a test vector field along $\G$. Then
$$
dF_a (u) = \frac{1}{2} \int_{\G} [(u\times \G') + (\G \times u')]\cdot a\ dt =  \int_{\G} \det (\G',a,u)\ dt = i_a \Omega(u),
$$
proving the first equality. The second equality follows from the identity $(a\times b)\cdot c = \det(a,b,c)$ for any triple of vectors.

For the last equality, let $v$ be a test vector field along $\G$. Then
$$
i_{\G\times a} \Omega (v) = \int_{\G} \det(\G',\G \times a,v)\ dt = 
\int_{\G} [(\G'\cdot a) (\G\cdot v) - (\G'\cdot \G) (a\cdot v)]  \ dt,
$$
where we use the identity 
$$
\det(a,b\times c,d) = (a\cdot c) (b\cdot d) - (a\cdot b) (c\cdot d)
$$
that holds for every quadruple of vectors. On the other hand, 
\begin{equation*}
\begin{split}
dJ_a(v) &= -  \int_{\G} [(v\cdot \G') (\G\cdot a) + (\G\cdot v') (\G\cdot a) + (\G\cdot \G') (v\cdot a)] \ dt\\
&= \int_{\G} [(\G'\cdot a) (\G\cdot v) - (\G'\cdot \G) (a\cdot v)]  \ dt,
\end{split}
\end{equation*}
where we used integration by parts. This completes the proof.
\proofend

Therefore, $a$ is a Hamiltonian vector field of the function $F_a$ with respect to $\Omega$. The constant vector field $a$ commutes with the filament flow and all higher flows $X_n$, and hence $F_a$, and therefore $A(\G)$, is an integral of the filament equation. For a discrete version of the filament equation, the area vector is a monodromy integral, see \cite{PSW}. 

\bigskip
{\bf Acknowledgments}. 
I am grateful to J. Langer and R. Perline for introducing me to the beautiful theory of the filament equation and related topics; many thanks to R. Perline for pointing out the connection between the bicycle model and the filament equation.
 It is a pleasure to acknowledge stimulating discussions on the subject of this note that I have had over the years with G. Bor, R. Foote,  T. Hoffmann, B. Khesin,  M. Levi,  V. Ovsienko, F. Pedit, U. Pinkall, C. Roger,   Yu. Suris, W. Thurston, E. Tsukerman,  A. Veselov, and V. Zharnitsky. The author was   supported by the NSF grants DMS-1105442 and  DMS-1510055.


\begin{thebibliography}{99}

\bibitem{AK} V. Arnold, B. Khesin. {\it Topological methods in hydrodynamics.} Springer-Verlag, New York, 1998.

\bibitem{Au} H. Auerbach. {\it Sur un probl\`eme de M. Ulam concernant l?\`equilibre des corps flottants}. Studia Math. {\bf 7} (1938), 121--142.

\bibitem{BMO1} J. Bracho, L. Montejano, D. Oliveros, {\it A classification theorem for Zindler carrousels.} J. Dynam. Control Systems {\bf 7} (2001),  367--384.

\bibitem{BMO2} J. Bracho, L. Montejano, D. Oliveros, {\it Carousels, Zindler curves and the floating body problem.} Period. Math. Hungar. {\bf 49} (2004),  9--23.

\bibitem{Ca} A. Calini,. {\it Recent developments in integrable curve dynamics.} Geometric approaches to differential equations, 56?99, Austral. Math. Soc. Lect. Ser., 15, Cambridge Univ. Press, Cambridge, 2000.

\bibitem{CI} A. Calini, T. Ivey. {\it B\"acklund transformations and knots of constant torsion.} J. Knot Theory Ramifications {\bf 7} (1998), 719--746.

\bibitem{Fi} D. Finn, {\it Which way did you say that bicycle went?}  Math. Mag. {\bf 77} (2004), 357--367.

\bibitem{FLT} R. Foote, M. Levi, S. Tabachnikov. {\it Tractrices, bicycle tire tracks, hatchet planimeters, and a 100-year-old conjecture}. Amer. Math. Monthly {\bf 120} (2013), 199--216.

\bibitem{Ho} T. Hoffmann. {\it Discrete Hashimoto surfaces and a doubly discrete smoke-ring flow}. Discrete differential geometry, 95--115, Oberwolfach Semin., 38, Birkh\"auser, Basel, 2008.

\bibitem{HPZ} S. Howe, Sean, M. Pancia, V. Zakharevich. {\it Isoperimetric inequalities for wave fronts and a generalization of Menzin's conjecture for bicycle monodromy on surfaces of constant curvature.} Adv. Geom. {\bf 11} (2011), 273--292.

\bibitem{LP} J. Langer, R. Perline. {\it Poisson geometry of the filament equation.} J. Nonlinear Sci. {\bf 1} (1991),   71--93. 

\bibitem{LP1} J. Langer, R. Perline. {\it The planar filament equation.} Mechanics day (Waterloo, ON, 1992), 171--180, Fields Inst. Commun., 7, Amer. Math. Soc., Providence, RI, 1996.

\bibitem{La} J. Langer. {\it Recursion in curve geometry.} New York J. Math. {\bf 5} (1999), 25--51.

\bibitem{LT} M. Levi, S. Tabachnikov, {\it On bicycle tire tracks geometry, hatchet planimeter, Menzin's conjecture, and oscillation of unicycle tracks}. Experiment. Math. {\bf 18} (2009),  173--186.

\bibitem{Le} M. Levi. {\it ``Bike tracks", quasi-magnetic forces, and the Schr\"odinger equation}. SIAM News {\bf 47},  June 2014.

\bibitem{MW} J. Marsden, A.  Weinstein.
{\it Coadjoint orbits, vortices, and Clebsch variables for incompressible fluids.} 
Phys. D {\bf 7} (1983),  305--323. 

\bibitem{PSW} U. Pinkall, B. Springborn, S. Weissmann, {\it A new doubly discrete analogue of smoke ring flow and the real time simulation of fluid flow}. J. Phys. A {\bf 40} (2007),  12563--12576.

\bibitem{RS} C. Rogers, W. Schief. {\it B\"acklund and Darboux transformations.} Cambridge University Press, Cambridge, 2002.

\bibitem{Ta} S. Tabachnikov, {\it Tire track geometry: variations on a theme}. Israel J. Math. {\bf 151} (2006), 1--28.

\bibitem{TT} S. Tabachnikov, E. Tsukerman. {\it On the discrete bicycle transformation.} Publ. Mat. Urug. {\bf 14} (2013), 201--219. 

\bibitem{We1} F. Wegner, {\it Floating bodies of equilibrium.} Stud. Appl. Math. {\bf 111} (2003),  167--183.

\bibitem{We2} F. Wegner, {\it Floating Bodies of Equilibrium I}. arXiv:physics/0203061.

\bibitem{We3} F. Wegner, {\it Floating Bodies of Equilibrium II}. arXiv:physics/0205059.

\bibitem{We4} F. Wegner, {\it Floating Bodies of Equilibrium. Explicit Solution}. arXiv:physics/0603160.

\bibitem{We5} F. Wegner, {\it Floating Bodies of Equilibrium in 2D, the Tire Track Problem and Electrons in a Parabolic Magnetic Field}. arXiv:physics/0701241.

\bibitem{We6} F. Wegner, {\it Three problems -- one solution}. \url{http://www.tphys.uni-heidelberg.de/~wegner/Fl2mvs/Movies.html#float}

\bibitem{Zi} K. Zindler. {\it \"Uber konvexe Gebilde II}. Monatsh. Math. Phys. {\bf 31} (1921), 25--57.


\end{thebibliography}
\end{document}